\documentclass[12pt,doublespacing]{article}

\usepackage{fullpage,setspace}
\usepackage{amssymb,amsmath}

\def\be{\begin{equation}} % useful definitions
\def\ee{\end{equation}}
\newcommand{\beq}{\begin{eqnarray}}
\newcommand{\eeq}{\end{eqnarray}}
\newcommand{\nbeq}{\begin{eqnarray*}}
\newcommand{\neeq}{\end{eqnarray*}}
\def\D{\displaystyle}

\begin{document}

\title{Variance estimators in critical branching processes with non-homogeneous immigration}
\author{I. Rahimov \and George P. Yanev \thanks{
First author: Department of Mathematics and Statistics, Zayed University, United Arab Emirates;  Second author: Mathematics Department, University of Texas - Pan American, USA.}
}
\date{\empty}
\maketitle

\begin{center}
{\small \textbf{Abstract}}

\end{center}

\noindent     The asymptotic normality of conditional least squares estimators for the offspring
variance in critical branching processes with non-homogeneous immigration is established, under moment assumptions on both reproduction and immigration. The proofs use martingale techniques and weak convergence results in Skorokhod spaces. \newline

\noindent \textit{Key Words:} branching processes, immigration, least squares estimators, offspring variance, Skorokhod space.\newline

\noindent \textit{MSC:}\ Primary 60J80;
Secondary 62F12, 60G99.

\section{Introduction}

Kimmel and Axelrod (2002), Pakes (2003), and Haccou et al. (2005) survey applications of branching stochastic models in genetics, molecular biology, and microbiology. Yakovlev and Yanev (2006) point out that {\it in vivo} cell kinetics requires stochastic modeling of renewing cell populations with non-homogeneous immigration. In this line, Hyrien and Yanev (2010) model renewing cell populations where the experimentally observable cells are supplemented by unobservable cells, for example stem cells. They analyze the population of terminally differentiated oligodendrocytes of the central nervous system and the population of leukaemia cells. In both cases the cell population expands through both division of existing (progenitor) cells and differentiation of stem cells. This dynamics belongs to branching processes with non-homogeneous immigration. The population's  viability (Jagers and Harding, 2009) is preserved by allowing the immigration distribution to vary in time, increasing to infinity on average. For semi-stochastic models where immigration depends on the state of the process we refer to Cairns (2009).

In our case of time-dependent immigration (Rahimov, 1995),
$\{X_{n,i}\}_{n,i\ge 1}$ and $\{ \xi_{n}\}_{n\ge 1}$   are two families of independent, non-negative, and
integer valued random variables on a probability space $(\Omega, \mathfrak{F}, P)$. We consider discrete-time branching processes with time-dependent immigration defined recursively by:
\begin{equation}  \label{def}
Z_n=\sum_{i=1}^{Z_{n-1}}X_{n,i}+\xi_{n}, \qquad n\geq 1; \quad Z_0=0,
\end{equation}
where the $\{X_{n,i}\}_{n,i\ge 1}$ have a common distribution for all $n$ and $i$, and
the sequences $\{X_{n,i}\}_{n,i\ge 1}$ and $\{\xi_{n}\}_{n\ge 1}$ are independent of each other. $X_{n,i}$
is the number of offspring of the $i$th  individual from
the $(n-1)$th generation and $\xi_{n}$ the number of immigrants (or invaders,  Haccou et al., 2005)
joining the population at the time of birth of the $n$th generation. $Z_n$ represents the $n$th population
size and the independence of $\{X_{n,i}\}_{n,i\ge 1}$ from $\{\xi _{n}\}_{n\ge 1}$ implies
the independence of reproduction and immigration. In
contrast to usual branching process models with immigration, we do not assume $\{\xi _{n}\}_{n\ge 1}$ to be
identically distributed, but that the immigration rate varies from generation to generation.

Our goal is to estimate the offspring variance $b^{2}:={\rm Var} X_{1,1}$ assuming that the immigration mean and variance are known. We  estimate $b^2$, based on observing a single trajectory $\{ Z_1, Z_2, \ldots, Z_n\}$ as $n\to \infty$. For reproduction, we assume
\be \label{criticality}
{\rm E}X_{1,1}=1\quad  \mbox{and} \quad {\rm E}X^4_{1,1}<\infty.
\ee
 For immigration,
let $\alpha_n  :=  {\rm E}\xi_{n}$, $\beta^2_n:={\rm Var} \xi_{n}$,
and $\gamma^4_n  :=  {\rm Var} (\xi_n-\alpha_n)^2$ be finite
for every  $n$ and varying regularly at infinity functions of $n$ defined  by:
\be  \label{2}
\alpha_n = n^{\alpha}L_{\alpha}(n), \quad  \beta^2_n=n^{\beta}L_{\beta}(n),\quad \mbox{and}\quad  \gamma^4_n  =  n^{\gamma}L_{\gamma}(n),
\ee
where $\alpha$, $\beta$, $\gamma$ are non-negative and  $ \lim_{x\to \infty}L_{(\cdot)}(cx)/L_{(\cdot)}(x)=1$ for any $c>0$. %$L_{\alpha}(n)$, $L_{\beta}(n)$, $L_{\gamma}(n)$ are
 The immigration mean increases to infinity:
\be \label{alpha_beta}
\lim_{n\to \infty}\alpha_n=\infty,
\ee
and the immigration moments, in addition to Eq.~(\ref{2}), satisfy:
\be  \label{gamma}
\mbox{(i)} \lim_{n\to \infty} \frac{\beta^2_n}{n\alpha_n}=0 \quad \mbox{or}\quad \mbox{(ii)}\lim_{n\to \infty}\frac{\beta^2_n}{n\alpha^2_n}=0\quad \mbox{and}\quad \lim_{n\to \infty} \frac{n\alpha_n\beta^2_n}{\gamma^4_n}=0.
\ee
Eq.~(\ref{gamma}) holds true if $\beta<\alpha+1$ or ($\beta<2\alpha+1$ and $\alpha+\beta+1<\gamma$). Finally, setting $\widetilde{\eta}_{n}:=(\xi_n-\alpha_n)^2-\beta^2_n$, we assume that $\{\widetilde{\eta}_n\}_{n\ge 1}$ satisfies the Lindeberg condition:
\be \label{moment_ass1}
\lim_{n\to \infty}\sum_{k=1}^n E\left(\widetilde{\eta}^2_k\chi (|\widetilde{\eta}_k|>\varepsilon)\right)=0,
\ee
where $\chi (\textit{A})$ denotes the indicator of the event $\textit{A}$.

To construct the conditional least squares estimator (CLSE) for the offspring variance, consider $\Im (k)$ the $\sigma$-algebra generated by $\{Z_j\}_{j=1}^k$. Denote ${\rm E}_k(\cdot):={\rm E}(\cdot |\Im (k))$ and define:
\be \label{M_def}
M_k   := Z_k-{\rm E}_{k-1}Z_k
 =  Z_k-Z_{k-1}-\alpha_k.
\ee
Denote $V_k := M^2_k-{\rm E}_{k-1}M^2_k$. The recurrence Eq.~(\ref{def}) yields the stochastic regression equation
\be \label{SRE}
M^2_k=b^2Z_{k-1}+\beta^2_k+V_k.
\ee
The error terms $V_k$ form
 a martingale difference sequence $\{V_k\}_{k\ge 1}$ with respect to the filtration $\{\Im (k)\}_{k\ge 1}$, that is, $V_k$ is $\Im (k-1)$-measurable and ${\rm E}_{k-1}V_k=0$. If $\alpha_k$ and $\beta^2_k$ are known and minimizing  with respect to $b^2$ the sum of squares
 \be
  \sum_{k=1}^n V^2_k=\sum_{k=1}^n\left(M^{2}_k-b^2Z_{k-1}-\beta^2_k\right)^2,
 \ee
 the CLSE for $b^2$ is:
\be \label{estimator}
\widehat{b^{2}_n} := \frac{\D \sum_{k=1}^{n}(M^2_k-\beta^2_k)Z_{k-1}}{\D \sum_{k=1}^{n}Z^2_{k-1}}
  =  \frac{\D \sum_{k=1}^{n}\left((Z_k-Z_{k-1}-\alpha_k)^2-\beta^2_k\right)Z_{k-1}}{\D \sum_{k=1}^{n}Z^2_{k-1}}.
\ee
For $\{Z^\ast_n\}_{n\ge 0}$, a simple branching process with homogeneous immigration, assuming that both offspring and immigration means are known, the CLSE for the offspring variance is given by
\be \label{homog}
\widehat{b^{ 2\ast}_n}:=\frac{\D \sum_{k=1}^{n}(M^\ast_k)^2(Z^\ast_{k-1}-\bar{Z}^\ast_n)}{\D \sum_{k=1}^{n}(Z^\ast_{k-1}-\bar{Z}^\ast_n)^{2}},
\ee
where
$\bar{Z}^\ast_n:=n^{-1}\sum_{k=1}^{n}Z^\ast_{k-1}$ and $M^\ast_k$ is defined with the first equality of Eq.~(\ref{M_def}), replacing $Z_k$ by $Z^\ast_k$.
 Yanev (1976/77) and Yanev and Tchoukova-Dantcheva (1986) pioneered the study of $\widehat{b^{2\ast}_n}$;
Winnicki (1991) proved limit theorems for $\widehat{b^{2\ast}_n}$
assuming finite fourth moments for both reproduction and homogeneous immigration (surveys in Dion, 1993 and  Yanev, 2008). Ma and Wang (2010) studied the case when these fourth moments may be infinite.

 We extend Winnicki's (1991) findings in the critical case ${\rm E}X_{1,1}=1$  allowing the immigration distribution to vary with time such that its mean $\alpha_n$ increases to infinity at a speed that could correspond to a near-critical branching population.  We prove a limit theorem, which completes the results of Rahimov (2008a) in that the limiting distribution of the CLSE for the offspring mean depends on the offspring variance $b^2$. We establish asymptotical normality of the CLSE defined by Eq.~(\ref{estimator}).

 We shall express $\widehat{b^{2}_n}-b^2$, suitably normalized, as a sum of martingale differences and functionals of $\{Z_k\}_{k=1}^n$. Then we study the asymptotic behavior of $\widehat{b^{2}_n}-b^2$ by applying either a limit theorem due to Rahimov (2007) giving the limit constant of some functionals of  $\{Z_k\}_{k=1}^n$ (Lemma~1) or a general limit theorem,  given as Lemma~2, for sequences of martingale differences in continuous time (Jacod and Shiryaev, 2003)  to each term of the sum.

Define $A_n:={\rm E}Z_n$,  $\tau^{2}_n:=\sum_{k=1}^{n}\gamma^4_k$, and
$\theta_n:= nA^2_n(nA_{n}^{2}+\tau_{n}^{2})^{-1}$. We have
\be \label{10+}
\lim_{n\to \infty}\theta_n=\lim_{n\to \infty}\frac{nA^2_n}{nA_{n}^{2}+\tau_{n}^{2}}=:\theta\in [0,1].
\ee
Equality in distribution is denoted by ``$\stackrel{d}{=}$".

\noindent{\bf Theorem 1}\ Assume Eq.~(\ref{criticality})-(\ref{moment_ass1}) hold true.
Then
\be \label{claim}
\lim_{n\to \infty}(\theta_nn)^{1/2}\left(\widehat{b^{2}_n}-b^2\right) \stackrel{d}{=} N(0,\sigma^2),
\ee
where $N(0,\sigma^2)$ is a normal random variable with zero mean and variance
\be \label{sigma}
\sigma^2=(2\alpha+3)^2\left(\theta\ \frac{2b^4}{4\alpha+5}+(1-\theta)\ \frac{\gamma+1}{2\alpha+3+\gamma}\right).
\ee
For the critical  process with homogeneous immigration, Winnicki  (1991)  established the weak limit of $\widehat{b^{2\ast}_n}$ with rate of convergence $n^{1/2}$. In Theorem 1 the convergence rate is $(\theta_n n)^{1/2}$, where $\lim_{n\to \infty}\theta_n=\theta\in [0,1]$ and  $\lim_{n\to \infty}\theta_nn=\infty$, provided $\lim_{n\to \infty}n^3\alpha^2_n/\gamma^4_n=0$.

\noindent{\bf Corollary 1}\ Under the assumptions of Theorem 1,

(i) If $\lim_{n\to \infty}n^2\alpha^2_n/\gamma^4_n=0$, then Eq.~(\ref{claim})-(\ref{sigma}) hold true with $\theta=0$.

(ii) If $\lim_{n\to \infty}n^2\alpha^2_n/\gamma^4_n=\infty$, then Eq.~(\ref{claim})-(\ref{sigma}) hold true with $\theta=1$.

\noindent{\bf Example 1 (Poisson immigration)}\ For $\{\xi_n\}_{n\ge 1}$ independent Poisson variables with mean
$\alpha_n=n^\alpha L_\alpha(n)=o(n)\to \infty$ as $n\to \infty$, Eq.~(\ref{alpha_beta})-(\ref{moment_ass1}) and the condition (ii) in  Corollary~1  satisfied, if Eq.~(\ref{criticality}) holds true, then Theorem~1 implies Eq.~(\ref{claim})-(\ref{sigma}) with $\theta=1$.

\noindent{\bf Example 2 (Neyman Type A immigration)}\ If Eq.~(\ref{criticality}) holds true and  for $\{\xi_n\}_{n\ge 1}$ independent with Neyman Type A distribution
given by
${\rm E}z^{\xi_n}=\exp\left(\lambda_n\left(e^{\varphi_n(z-1)}-1\right)\right)$, $|z|<1$,
if
$\lambda_n=n^\lambda L_\lambda(n)\to \infty$, $(\lambda\ge 0)$ and $\varphi_n=n^\varphi L_\varphi(n)$, $(\varphi\ge 0)$, then
(Johnson et al., 1993: 371)  the $r$th factorial moment $\mu_n(r)$ satisfies
$\mu_n(r)\sim \lambda^r_n\varphi^r_n$ for  $r\ge 1$ and $n\to \infty$.  Also $\{\widetilde{\eta}_n\}_{n\ge 1}$ satisfies Eq. (\ref{moment_ass1}).
If
 $0<\lambda+\varphi \le 1/2$,
 then Eq.~(\ref{alpha_beta}) and (\ref{gamma})(i) hold true and hence Eq.~(\ref{claim})-(\ref{sigma}) with $\theta=1$ and $\alpha=\lambda+\varphi$. If $1<\lambda+\varphi<3/2$, then  Eq.~(\ref{alpha_beta}) and  (\ref{gamma})(ii) are satisfied, which yields Eq.~(\ref{claim})-(\ref{sigma}) with $\theta=0$ and $\alpha=\lambda+\varphi$.

\section{Preliminaries}

The proof of the theorem uses auxiliary results
given in this section. ``$D$ " denotes convergence or equality in the Skorokhod space $D(\mathbb{R}_{+},
\mathbb{R})$, ``$P$ " probability, and ``$d$ " distribution.
The first lemma summarizes limit results for functionals of Eq.~(\ref{def}). Its proof is similar to that of Corollary 2 in Rahimov (2008b) and is omitted here.

{\bf Lemma 1}  For the critical process (\ref{def}) with $b^2<\infty$,
$\lim_{n\to \infty}\alpha_n=\infty$ and \newline $\lim_{n\to \infty}\beta^2_n(n\alpha^{2}_n)^{-1}=~0$, for
any continuous function $\Phi$ on  $\mathbb{R}_+$ and
any sequence $\{c_n\}_{n\ge 0}$, varying regularly at infinity with exponent $\rho\ge 0$,  for $t>0$ we have:
\be
\lim_{n\to \infty}\frac{1}{nc_n}\sum_{k=0}^{[nt]}c_k\Phi\left(\frac{Z_k}{A_n}\right)\stackrel{P}{=}\int_{0}^{t}u^{\rho}\Phi(u^{\alpha+1})\, du.
\ee

A necessary and sufficient condition for weak convergence in a Skorokhod space of a sequence of martingale differences (Jacod and Shiryaev, 2003: Theorem VIII.2.29; Isp$\grave{\mbox{a}}$ny et al., 2006) is:

{\bf Lemma 2 (CLT for martingales)}   For a sequence of martingale differences  $\{U_k^n\}_{k\ge 1}$, $n\ge 1$, with respect to a filtration $\{\Im_k^n\}_{k\ge 1}$, such that for all $\varepsilon>0$ and $t\ge 0$ the Lindeberg  condition
\be \label{Lin_C}
\lim_{n\to \infty}\sum_{k=1}^{[nt]}{\rm E}\left((U_k^n)^2\chi (|U_k^n|>\varepsilon)\ |\ \Im_{k-1}^n\right)\stackrel{P}{=}0
\ee
holds true. Then
\be
\lim_{n\to \infty}\sum_{k=1}^{[nt]}U_k^n\stackrel{D}{=}U(t),
\ee
where $U(t)$ is a continuous Gaussian martingale with mean zero and covariance function $C(t)$ if and only if for every $t\ge 0$,
\be \label{14}
\lim_{n\to \infty}\sum_{k=1}^{[nt]}{\rm E}\left((U_k^n)^2\ |\ \Im_{k-1}^n\right)\stackrel{P}{=}C(t).
\ee
We use a tilde to indicate that a random variable $\zeta$ is centered around its mean, that is  $\widetilde{\zeta}=\zeta-{\rm E}\zeta$. Denote  $Y_{k,i}:=\widetilde{X}_{k,i}^2$, $S_{j}:=\sum_{i=1}^{j}\widetilde{X}_{k,i}$, and  recall that $\eta_n:=(\xi_n-\alpha_n)^2$. The  expansion of the error term  $V_k$ plays a key role in our analysis:
\beq \label{V_k_repres}
V_k  & = & M^2_k-{\rm E}_{k-1}M^2_k \\
& = &
(Z_k-Z_{k-1}-\alpha_k)^{2}-b^{2}Z_{k-1}-\beta^2_k \nonumber \\
 & = & \left(2\sum_{j=2}^{Z_{k-1}}\widetilde{X}_{k,j}S_{j-1}+\widetilde{\eta}_{k}\right)
+2\widetilde{\xi}_{k}\sum_{i=1}^{Z_{k-1}}\widetilde{X}_{k,i}+\sum_{i=1}^{Z_{k-1}}\widetilde{Y}_{k,i}
\nonumber \\
    & =: &
    V_k^{(1)}+V_k^{(2)}+V_k^{(3)}.\nonumber
\eeq
Denoting
$
V^{(i)}_{n}(t)
     :=  H^{-1}_n\sum_{k=1}^{[nt]}V^{(i)}_{k}$, $1\le i\le 3$, where $H_{n}^{2}=nA_{n}^{2}+\tau_{n}^{2}$,
we decompose the normalized sum of the error terms into three parts:
\be \label{main_decomp}
V_n(t):=\frac{1}{H_n}\sum_{k=1}^{[nt]}V_{k}=V^{(1)}_{n}(t)+ V^{(2)}_{n}(t)+ V^{(3)}_{n}(t), \quad t>0.
\ee
We shall show that the asymptotic behaviour of $V_n(t)$ as $n\to \infty$ is governed by $V_n^{(1)}(t)$, while the contributions of $V_n^{(2)}(t)$ and $V_n^{(3)}(t)$ are negligible. Define
\be \label{process_V}
V(t):= W\left(\theta\frac{2b^4}{2\alpha+3}t^{2\alpha+3}+(1-\theta)t^{\gamma+1}\right),\quad t>0,
\ee
where $W(t)$ is a standard Wiener process and $\theta$ is the limiting constant in Eq.~(\ref{10+}).

{\bf Proposition 1.} If Eq.~(\ref{criticality})-(\ref{gamma}), then for every $t>0$
\be \label{prop1}
\lim_{n\to \infty}V^{(1)}_{n}(t)
\stackrel{D}{=}V(t).
\ee
{\bf Proof.}\ Denote $T_k:=2\sum_{j=2}^{Z_{k-1}}\widetilde{X}_{k,j}S_{j-1}$. The independence of $\{X_{n,i}\}_{n,i\ge 1}$ and Lemma~1 yield
\beq \label{new_V}
\sum_{k=1}^{[nt]}{\rm E}_{k-1}\left(\frac{V_k^{(1)}}{H_n}\right)^{2}
%& = &
%\frac{1}{H^2_n}\sum_{k=1}^{[nt]}{\rm E}_{k-1}\left(T_k+\widetilde{\eta}_k\right)^2 \\
    &  = &
\frac{1}{H^2_n}\sum_{k=1}^{[nt]}{\rm E}_{k-1}T^2_k
+ \frac{1}{H^2_n}\sum_{k=1}^{[nt]}{\rm E}_{k-1}\widetilde{\eta}_k^2  \\
    & = &
     \frac{2b^4\theta_n}{nA^2_n}\sum_{k=1}^{[nt]} Z_{k-1}(Z_{k-1}-1)+\frac{1-\theta_n}{\tau_n^2}\sum_{k=1}^{[nt]}\gamma_k^4
      \nonumber\\
   & \stackrel{P}{\rightarrow} &
\frac{2b^4\theta}{2\alpha+3}t^{2\alpha+3}+  (1-\theta)t^{\gamma+1}. \nonumber
\eeq
Then $V_k^{(1)}/H_n$ satisfies Eq.~(\ref{14}) with $C(t)=2b^4\theta t^{2\alpha+3}(2\alpha+3)^{-1}+  (1-\theta)t^{\gamma+1}$.
We shall verify Eq.~(\ref{Lin_C}). Indeed,
\beq
\lefteqn{\hspace{-2cm}\sum_{k=1}^{[nt]}{\rm E}_{k-1}\left(\left(\frac{V_k^{(1)}}{H_n}\right)^2 \chi\left( \left|\frac{V_k^{(1)}}{H_n}\right|>\varepsilon\right)\right)
 \le
\frac{1}{H^2_n}\sum_{k=1}^{[nt]}{\rm E}_{k-1}\left((T^2_k + \widetilde{\eta}^2_k)\chi\left( \left|T_k+\widetilde{\eta}_k\right|>\varepsilon H_n\right)\right) } \nonumber \\
& \le &
\frac{1}{H^2_n}\sum_{k=1}^{[nt]}{\rm E}_{k-1}\left(T^2_k \ \chi\left( \left|\widetilde{\eta}_k\right|>\frac{\varepsilon H_n}{2}\right)+T^2_k \ \chi\left( \left|T_k\right|>\frac{\varepsilon H_n}{2}\right)\right) \nonumber \\
& & +
\frac{1}{H^2_n}\sum_{k=1}^{[nt]}{\rm E}_{k-1}\left(\widetilde{\eta}^2_k  \ \chi\left( \left|\widetilde{\eta}_k\right|>\frac{\varepsilon H_n}{2}\right)+\widetilde{\eta}^2_k \ \chi\left( \left|T_k\right|>\frac{\varepsilon H_n}{2}\right)\right)\nonumber\\
& =: & I_1(n)+I_2(n)+I_3(n)+I_4(n).
\eeq
The independence assumption, Chebyshev inequality, and Lemma 1 imply
\beq \label{Ch_1}
I_1(n) +I_4(n)& = & \frac{1}{H^2_n}\sum_{k=1}^{[nt]}\left({\rm E}_{k-1}(T^2_k)  P\left( \left|\widetilde{\eta}_k\right|>\frac{\varepsilon H_n}{2}\right)+
{\rm E}_{k-1}(\widetilde{\eta}^2_k) \ P\left( \left|T_k\right|>\frac{\varepsilon H_n}{2}\right)\right)\nonumber\\
    & \le &
\frac{4}{\varepsilon^2 H^4_n}\sum_{k=1}^{[nt]}\left({\rm Var}_{k-1}T_k \ {\rm Var}\widetilde{\eta}_k+{\rm Var}\widetilde{\eta}_k \ {\rm Var}_{k-1}T_k\right)  \nonumber \\
& = &
\frac{8b^4\theta_n(1-\theta_n)}{\varepsilon^2 \tau^2_n}\frac{1}{n}\sum_{k=1}^{[nt]}\frac{Z_{k-1}(Z_{k-1}-1)}{A^2_n}\gamma^4_k \nonumber \\
& \stackrel{P}{\rightarrow} & 0.
\eeq
For $I_3(n)$, referring to the Lindeberg condition in Eq.~(\ref{moment_ass1}), we have
\beq
I_3(n) & = & \frac{1}{H^2_n}\sum_{k=1}^{[nt]}{\rm E}_{k-1}\left(\widetilde{\eta}^2_k  \ \chi\left( \left|\widetilde{\eta}_k\right|>\frac{\varepsilon H_n}{2}\right)\right)\\
        & \le &
 \frac{1-\theta_n}{\tau^2_n}\sum_{k=1}^{[nt]}{\rm E}_{k-1}\left(\widetilde{\eta}^2_k  \ \chi\left( \left|\widetilde{\eta}_k\right|> \frac{\varepsilon\tau_n}{2}\right)\right)\nonumber \\
 & \stackrel{P}{\rightarrow} & 0. \nonumber
\eeq
It remains to show $\lim_{n\to \infty}I_2(n)\stackrel{P}{=}  0$.
From Burkholder inequality, for $\delta>0$ and $p>1$,
\be \label{Burk}
{\rm E}\left|\sum_{j=2}^{r}\widetilde{X}_{k, j}S_{j-1}\right|^{2p}
\leq
D_{1,p}{\rm E}\left|\sum_{j=2}^{r}\widetilde{X}_{k,j}^2S_{j-1}^{2}\right|^{p},
\ee
where $D_{1,p}>0$ depends  on $p$ only. Using
Minkovski inequality,
\be \label{I_1}
\left({\rm E}\left|\sum_{j=2}^{r}\widetilde{X}_{k,j}^2S_{j-1}^{2}\right|^{p}\right)^{\frac{1}{p}}
     \le
    \sum_{j=2}^{r}\left({\rm E}\left(|\widetilde{X}_{k, j}|^{2p}|S_{j-1}|^{2p}\right)\right)^{\frac{1}{p}}
     =
    b^{\frac{1}{p}}_{2p}\sum_{j=2}^{r}\left({\rm E}|S_{j-1}|^{2p}\right)^{\frac{1}{p}},
\ee
where,  by assumption, $b_{2p}:={\rm E}|\widetilde{X}_{1,1}|^{2p}<\infty$ for $0<p\le 2$. Similarly, we obtain
\be \label{I_2}
{\rm E}|S_{j-1}|^{2p}  \le
D_{2,p}{\rm E}\left|\sum_{i=1}^{j-1}\widetilde{X}_{k,i}^2\right|^{p}
   \le   D_{2,p}b_{2p}(j-1)^p,
\ee
where $D_{2,p}>0$ depends on
$p$ only.
From Eq.~(\ref{Burk})-(\ref{I_2}),
\be \label{25}
{\rm E}\left|\sum_{j=2}^{r}\widetilde{X}_{k, j}S_{j-1}\right|^{2p}
 \leq
D_{1,p}D_{2,p}b_{2p}\left[\sum_{j=2}^{r}(j-1)\right]^{p}
 \leq
D_{p}b_{2p}r^{2p},
\ee
where $D_{p}>0$ depends on
$p$ only. Applying Eq.~(\ref{25}) with $p=(2+\delta)/2$ to ${\rm E}_{k-1}|T_k|^{2+\delta}$, for every $\varepsilon>0$ and $\delta\ge 0$:
\beq \label{LC}
I_2(n) & \le &
\frac{1}{\varepsilon^\delta H_n^{2+\delta}}\sum_{k=1}^{[nt]}{\rm E}_{k-1}|T_k|^{2+\delta}
\\
& \leq &
\frac{D_{\delta}b_{2+\delta}}{\varepsilon^{\delta}H^{2+\delta}_n}\sum_{k=1}^{[nt]}Z^{2+\delta}_{k-1}\nonumber \\
& = &
\frac{D_{\delta}b_{2+\delta}\theta_n}{\varepsilon^{\delta}}\left(\frac{A_n}{H_n}\right)^\delta \frac{1}{n}
\sum_{k=1}^{[nt]}\left(\frac{Z_{k-1}}{A_n}\right)^{2+\delta} \nonumber \\
& \stackrel{P}{\to} & 0, \nonumber
    \eeq
    where we have used Lemma 1 and $\lim_{n\to \infty} A_n/H_n=0$.
 Hence, $V_k^{(1)}/H_n$ satisfies Eq.~(\ref{Lin_C}) and because   $V^{(1)}_k$ is a martingale difference with respect to $\{\Im_k\}_{k\ge 1}$, all assumptions of Lemma~2 are verified. Lemma 2 implies Eq.~(\ref{prop1}).

 We derive the limit of the normalized sum
in Eq.~(\ref{main_decomp}).

{\bf Proposition 2}\  If Eq.~(\ref{criticality})-(\ref{moment_ass1}), then for every $t>0$
\be \label{prop2new}
\lim_{n\to \infty}V_n(t)
\stackrel{D}{=}  V(t).
\ee
{\bf Proof.}
From Eq.~(\ref{main_decomp}) and Proposition 1, in order to prove Eq.~(\ref{prop2new}), it remains to establish that  both  $V^{(2)}_{n}(t)$ and $V^{(3)}_{n}(t)$ are asymptotically negligible.
First we show
\be \label{V_4_P}
\lim_{n\to \infty}V^{(3)}_{n}(t)
\stackrel{P}{=}0.
\ee
%Let us verify the assumptions Eq.~(\ref{Lin_C}) and Eq.~(\ref{14}) for $U^n_k=V^{(3)}_{k}/H_n$.
Observing  that $\lim_{n\to \infty}nA_n/H^2_n= 0$ and applying Lemma 1, we obtain
\beq \label{24}
\sum_{k=1}^{[nt]}{\rm E}_{k-1}\left(\left(\frac{V_k^{(3)}}{H_n}\right)^2 \chi\left( \left|\frac{V_k^{(3)}}{H_n}\right|>\varepsilon\right)\right)
 & \le &  \frac{1}{H^2_n}\sum_{k=1}^{[nt]}{\rm E}_{k-1}\left(V_{k}^{(3)}\right)^2   \\
    & = &
    \frac{\theta_n}{A_n}\frac{{\rm E}\widetilde{Y}_{1,1}^2}{n}\sum_{k=1}^{[nt]}\frac{Z_{k-1}}{A_n}\nonumber \\
    & \stackrel{P}{\to} & 0, \nonumber
\eeq
where ${\rm E}\widetilde{Y}_{1,1}^2={\rm E}X^4_{1,1}-4{\rm E}X^3_{1,1}-4b^4+3b^2+3<\infty$ (due to Eq.~(\ref{criticality})).
Then, Eq.~(\ref{Lin_C}) and Eq.~(\ref{14}) with $C(t)\equiv 0$ hold true.
Lemma~2 yields $V^{(3)}_{n}(t)
\stackrel{D}{\to}0$ and hence Eq.~(\ref{V_4_P}).

As for $V^{(2)}_{n}(t)$,
applying Lemma 1, we obtain
\beq
\sum_{k=1}^{[nt]}{\rm E}_{k-1}\left(\frac{V_k^{(2)}}{H_n}\right)^2 \chi\left( \left|\frac{V_k^{(2)}}{H_n}\right|>\varepsilon\right)
  & \le &  \frac{1}{H^2_n}\sum_{k=1}^{[nt]}{\rm E}_{k-1}\left(V_{k}^{(2)}\right)^2  \\
    & = &
 \frac{\beta^2_n\theta_n}{A_n}\frac{4b^2}{n\beta^2_n}\sum_{k=1}^{[nt]}\beta^2_k\frac{Z_{k-1}}{A_n}\nonumber \\
 & \stackrel{P}{\to} & 0, \nonumber
    \eeq
where, according to Eq.~(\ref{gamma}), $\lim_{n\to \infty}\beta^2_n\theta_n/A_n=0$.
Therefore, Eq.~(\ref{Lin_C}) and (\ref{14}) with $C(t)\equiv~ 0$ are satisfied.
From Lemma~2 again, we have
$\lim_{n\to \infty}V^{(2)}_{n}(t)\stackrel{D}{=}0$ and hence
\be \label{V_2_P}
\lim_{n\to \infty}V^{(2)}_{n}(t)\stackrel{P}{=}0.
\ee
Eq.~(\ref{prop2new}) follows from Proposition 1, Eq.~(\ref{V_4_P}), Eq.~(\ref{V_2_P}),
and Slutsky theorem.

\vspace{-0.5cm}\section{Proof of Theorem 1}
Recalling Eq.~(\ref{estimator}) and the definition of $V_k$,  we have
\be \label{4}
\widehat{b^{2}_n}-b^{2}  =
\frac{\sum_{k=1}^{n}\left(M^{2}_k-\beta^2_k-b^2Z_{k-1}\right)Z_{k-1}}
{\sum_{k=1}^{n}Z^2_{k-1}}
     =
    \frac{\sum_{k=1}^{n}V_kZ_{k-1}}
{\sum_{k=1}^{n}Z^2_{k-1}}.
\ee
First, we examine the asymptotic behavior as $n\to \infty$ of the numerator in Eq.~(\ref{4}).
Theorem~3.1 in Rahimov~(2009) and Proposition~2 for $t>0$ imply:
\be
  \lim_{n\to \infty}Z_{n}(t)=\lim_{n\to \infty}\frac{Z_{[nt]}}{A_n}\stackrel{D}{=}t^{\alpha+1}\quad \mbox{and}\quad \lim_{n\to \infty} V_{n}(t)=\lim_{n\to \infty}\frac{1}{H_{n}}\sum_{k=1}^{[nt]}V_{k}\stackrel{D}{=}
    V(t),
\ee
where the convergence is  on $D(\mathbb{R}_{+},
\mathbb{R}_{+})$ and $D(\mathbb{R}_{+},
\mathbb{R})$, respectively.  Both limiting processes are
continuous, referring to
Theorem 2.2 in Kurtz and Protter (1991), so that
\begin{equation}\label{3.4}
    \lim_{n\to \infty}\left(Z_{n}(t),\;V_{n}(t),\;\int_{0}^{1}Z_{n}(u)\, dV_{n}(u)\right)
    \stackrel{D}{=}\left(t^{\alpha+1},\;V(t),\;\int_{0}^{1}u^{\alpha+1}\, dV(u)\right),
\end{equation}
on
$D(\mathbb{R}_{+}, \mathbb{R}_{+}\times\mathbb{R}\times\mathbb{R})$.
Eq.~(\ref{3.4}) and the continuous mapping theorem (Billingsley, 1968: Theorem 5.5) yield:
\beq \label{Prop4}
   \frac{1}{H_nA_n}\sum_{k=1}^{n}V_kZ_{k-1}
   & = & \sum_{k=1}^{n-1}Z_n\left(\frac{k}{n}\right)\left(V_n\left(\frac{k+1}{n}\right)-V_n\left(\frac{k}{n}\right)\right)
        \\
   & = & \int_{0}^{1}Z_{n}(u)\, dV_{n}(u)\nonumber \\
   & \stackrel{d}{\rightarrow} & \int_{0}^{1}u^{\alpha+1}\, dV(u).\nonumber
\eeq
For the denominator in  Eq.~(\ref{4}), applying Lemma 1 with $t=1$, we have
\be \label{C_3}
\lim_{n\to \infty} \frac{1}{nA^2_n}\sum_{k=1}^n Z^2_{k-1}
    \stackrel{P}{=}
    \frac{1}{2\alpha+3}.
    \ee
Finally, Eq.~(\ref{Prop4}) and (\ref{C_3}) and Slutsky theorem yield:
\be \label{Slutsky}
\lim_{n\to \infty} \frac{nA_n}{H_n}\left(\widehat{b^{2}_n}-b^2\right)
  =   \lim_{n\to \infty}\frac{ (H_nA_n)^{-1}\sum_{k=1}^n V_kZ_{k-1}}
        { (nA^2_n)^{-1}\sum_{k=1}^nZ^2_{k-1}}
 \stackrel{d}{=}
    (2\alpha+3)\int_{0}^{1}u^{\alpha+1}\, dV(u),
    \ee
which implies the normality of the limit in Eq.~(\ref{claim}).
Ito's formula yields:
\be  \label{Ito}
\int_0^1u^{\alpha+1}\, dV(u)  =  V(1)-(\alpha+1)\int_0^1u^{\alpha} V(u)\, du
    =  (\alpha+1)\int_0^1 (V(1)-V(u))u^\alpha \, du.
\ee
The variable  $\zeta:=\int_0^1 (V(1)-V(u))u^\alpha \, du$ is a normally distributed with zero mean and
\be
{\rm E}\zeta^2=\int_0^1\int_0^1s^\alpha t^\alpha {\rm E}((V(1)-V(t))(V(1)-V(s))) \, ds\, dt.
\ee
For $0\le s\le t\le 1$,
\beq
\lefteqn{{\rm E}((V(1)-V(t))(W(1)-W(s)))}\\
    & = &
    {\rm E}((V(1)-V(t))^2+(V(1)-V(t))(V(t)-V(s)))\nonumber \\
    & = & \theta\frac{2b^4}{2\alpha+3}\left(1-t^{2\alpha+3}\right)+(1-\theta)\left(1-t^{\gamma+1}\right), \nonumber
\eeq
we calculate
\be \label{var_zeta}
{\rm E}\zeta^2= \frac{1}{(\alpha+1)^2}\left( \theta\frac{2b^4}{4\alpha+5}+(1-\theta)\frac{\gamma+1}{2\alpha+\gamma+3}\right),
\ee
which, taking into account Eq.~(\ref{Slutsky})-(\ref{Ito}), implies Eq.~(\ref{sigma}) and completes the proof.

\section{Conclusion} Branching processes with time-dependent immigration are encountered in a variety of applications
to population biology. We study conditional least-squares estimators for the offspring variance in the critical case, assuming
that the immigration mean increases to infinity over time. Theorem 1 establishes the asymptotic normality of the proposed estimators
with convergence rate $(\theta_n n)^{1/2}$, where $\lim_{n\to \infty}\theta_n=\theta\in [0,1]$. A next question
concerns the conditional consistency of the estimators with different weights.

\vspace{0.3cm} {\bf ACKNOWLEDGEMENTS} We thank both referees for their helpful suggestions.

\vspace{0.3cm} {\bf REFERENCES}

 \noindent  Billingsley, P. (1968). {\it Convergence of Probability Measures.} New York: John Wiley.

\noindent Cairns, B.J. (2009).  Evaluating the Expected Time to Population
Extinction with Semi-

Stochastic Models. {\it Mathematical Population Studies, 16}: 199-220.

\noindent Dion, J.-P. (1993). Statistical inference for discrete time branching processes. In A. Obretenov,

V. Stefanov (Eds.), {\it Lecture Notes of the Seventh International Summer School on

Probability Theory and Mathematical Statistics.} Singapore: Sci. Culture

Technology Publishing, pp. 60-121.

\noindent  Haccou, P., Jagers, P., and Vatutin, V.A. (2005). {\it Branching Processes: Variation, Growth,

and Extinction of Populations.} Cambridge: Cambridge University Press.

\noindent  Hyrien, O. and Yanev, N.M. (2010). Modelling cell kinetics using branching processes with

nonhomogeneous Poisson immigration. {\it Proceedings of the Bulgarian Academy of
Sci-

ences, 63}(10): 1405-1414.

\noindent  Isp$\grave{\mbox{a}}$ny M., Pap G., and Van Zuijlen M.C.A. (2006). Critical branching mechanisms with

immigration and Ornstein-Uhlenbeck type diffusions. {\it Acta Scientiarum Mathematica-

rum, 71}: 821-850.

\noindent  Jacod, J. and  Shiryaev, A.N. (2003). {\it Limit Theorems for Stochastic Processes.} Berlin:

Springer.

\noindent Jagers, P. and Harding, K.C. (2009). Viability of Small Populations Experiencing
Recurring

Catastrophes. {\it  Mathematical Population Studies, 16}: 177-198.

\noindent  Johnson, N.L., Kotz, S., and Kemp, A.W. (1993). {\it Univariate Discrete Distributions (2nd

Edition).} New York: John Wiley.

\noindent  Kimmel, M. and Axelrod, D.E. (2002). {\it Branching Processes in Biology.} New York: Springer.

\noindent  Kurtz T. G. and Protter P. (1991). Weak limit theorems for stochastic integrals and stochas-

tic differential equations. {\it Annals of Probability, 19}(3): 1035-1070.

\noindent  Ma, Ch. and Wang, L. (2010). On estimation of the variances for critical branching processes

with immigration. {\it Journal of  Applied Probability, 47}: 526-542.

\noindent  Pakes, A.G. (2003). Biological applications of branching processes. In D.N. Shanbhag, C.R.

Rao (Eds.), {\it Handbook of Statistics 21, Stochastic Processes: Modeling and Simulation.}

Amsterdam: Elsevier, pp. 693-773.

\noindent  Rahimov, I. (1995). {\it Random Sums and Branching Stochastic Processes. Lecture Notes in

Statistics 96.} New York: Springer.

\noindent  Rahimov, I. (2007). Functional limit theorems for critical processes with immigration.

{\it Advances in Applied Probability, 39}: 1054-1069.

\noindent  Rahimov, I. (2008a). Asymptotic distribution of the CLSE in a critical process with immi-

gration. {\it Stochastic Processes and  Applications, 118}: 1892-1908.

\noindent   Rahimov, I. (2008b). Deterministic approximation of a sequence of nearly critical branching

processes. {\it Stochastic Analysis and  Applications, 26}: 1013-1024.

\noindent  Rahimov, I. (2009). Asymptotially normal estimators for the offspring mean in the branching

process with immigration. {\it Communications in Statistics - Theory and Methods, 38}:13-28.

\noindent  Winnicki, J. (1991). Estimation of the variances in the branching process with immigration.

{\it Probability Theory and Related Fields, 88}: 77-106.

\noindent  Yakovlev, A. and Yanev, N.M. (2006). Branching stochastic processes with immigration in

analysis of renewing cell populations. {\it Mathematical Biosciences, 203}: 37-63.

\noindent  Yanev, N.M. (1976/77). Estimates of variances in a subcritical branching process with

immigration (in Russian, English summary). {\it Annuaire de L'Universit$\acute{e}$ de Sofia Kliment

Ohridski, Facult$\acute{e}$ de math$\acute{e}$matiques et m$\acute{e}$canique, 71}(2): 39-44.

\noindent  Yanev, N.M. (2008). Statistical inference for branching processes. In M. Ahsanullah and G.

Yanev (Eds.), {\it Records and Branching Processes.} New York: Nova Science Publishers,

pp. 147-172.

\noindent  Yanev, N.M. and Tchoukova-Dantcheva, S. (1986). Limit theorems for estimators of vari-

ances in a branching process with immigration. {\it Serdica Mathematical Journal, 12}(2):

143-153 (in Russian).

\end{document}